\documentclass{amsart}
\usepackage{amsmath,amssymb,amsthm,latexsym,graphicx,url}

\theoremstyle{plain}
\newtheorem{lemma}{Lemma}[section]
\newtheorem{proposition}[lemma]{Proposition}
\newtheorem{theorem}[lemma]{Theorem}
\theoremstyle{definition}

\newtheorem{example}[lemma]{Example}

\theoremstyle{remark}

\def\A{\Bbb A}
\def\P{\Bbb P}
\def\Q{\Bbb Q}
\def\dim{{\rm dim}}

\title{Plain Varieties}
\thanks{
The authors have been supported by the Austrian Science Fund (FWF)
in the frame of project P18992 and project SFB 1303, and by the
Austrian Exchange Service (\"OAD) in the frame of ``Acciones Integradas''.}

\author{G\'abor Bodn\'ar}
\address{RISC Software Ges.m.b.H, Hagenberg, Austria, Gabor.Bodnar@risc.uni-linz.ac.at}
\author{Herwig Hauser}
\address{Universit\"at Vienna and Innsbruck, Austria, herwig.hauser@univie.ac.at}
\author{Josef Schicho}
\address{Johann Radon Institut, \"Osterr. Akad. Wiss., Linz, Austria, josef.schicho@oeaw.ac.at}
\author{Orlando Villamayor}
\address{Universidad Aut\'onoma, Madrid, Spain, villamayor@uam.es}

\begin{document}

\maketitle

\begin{abstract}
\noindent 
Algebraic varieties which are locally isomorphic 
to open subsets of affine space will be called {\em plain}.  
Plain varieties are smooth and rational. The converse is true for curves 
and surfaces, and unknown in general. It is shown that
plain varieties are stable under blowup in smooth centers.
\\


\end{abstract}


\section{Introduction}
\label{sec:intro}

Differential manifolds are obtained by gluing open subsets
of Euclidean space through differentiable isomorphisms. 
In the algebraic category, the situation is different: 
there are smooth algebraic varieties that cannot be covered 
by open subsets of affine space. 
Any non-rational variety provides an example. 

It is nevertheless always possible to cover a smooth variety  by open subsets which are isomorphic to smooth hypersurfaces 
of affine space. This is a generalization of the theorem of the primitive element (cf.÷ Proposition~\ref{thm:emb}). The isomorphisms can be
obtained by generic projections onto linear subspaces of the
appropriate di\-mension.

Those varieties which are locally isomorphic to {\em open} subsets 
of $\Bbb A^n$ will be called {\em plain}. They form a very natural class 
of varieties for which local computations are particularly efficient. 
The affine charts are not just any  smooth subvarieties of affine space 
but open, hence dense subsets. 
This allows one to work with affine coordinates. Moreover, the affine 
coordinate rings have unique factorization. The concept of plain varieties 
also appears in \cite{Akbulut:92} as ``algebraic spaces'' and in 
\cite{Manin:74} as ``special varieties'', where they
are used for the definition of $R$-equivalence.

First examples of plain varieties are graphs of polynomial maps 
$\A^n\to\A^m$, projective spaces, Grassmannians,
and their cartesian products -- the class of plain varieties is obviously
closed under finite products. More generally, it is easily seen that 
smooth toric varieties are plain (see section~\ref{sec:plain}).

We prove in Theorem~\ref{thm:main} that the class of 
plain varieties is closed under blowup in smooth centers, 
provided that the base field is infinite. This is not immediate 
if the global structure of the center is complicated, e.g., if the center has large genus. By Proposition~\ref{thm:emb} we may embed the center locally as a hypersurface in an affine subspace and then use the explicit ring-theoretic description of blowups by subrings of the function field of the variety.

An important application is the following: In characteristic zero, 
it is possible to prove the embedded resolution of singularities
by blowing up the ambient space along smooth centers in the strict transform of the variety. In this way we
can compute embedded resolutions within plain ambient spaces. As these can be covered by open subsets of affine space, the number of ambient variables does not increase when passing to the transform of the variety.
For instance, we can implement Villamayor's algorithm \cite{Villamayor:89} 
or the algorithm of Bierstone and Milman \cite{Bierstone_Milman:97} for
the reso\-lution of a singular hypersurface in such a way that the final 
(and all intermediate) results are covered by {\em hypersurfaces}
in the literal sense, each of them given by one equation in 
independent variables. Indeed, this algorithm proceeds by blowing
up the ambient space along smooth centers; and if the given ambient space
is plain, then all subsequent ambient spaces are plain, too.

Giving necessary and sufficient criteria for a variety to be plain
seems to be a much harder question. Obviously, plain varieties are rational
and smooth. 
We do not know of an example of a smooth rational variety which
is not plain.


\section{Smooth varieties are locally hypersurfaces}
\label{sec:emb}

Let $k$ be a field, which is assumed to be infinite unless otherwise specified.
We work in the category of varieties over $k$,
as defined in \cite{Hartshorne:77}:
a variety is an integral scheme over $k$ of finite type (in particular, reduced and irreducible),
and a morphism is a regular map.
Affine space $\A^n_k$ is the spectrum of a polynomial ring $R$ 
in $n$ algebraically independent variables over $k$. 
Any generator system $x_1, \ldots ,x_n$ of $R$ will be called 
a {\em coordinate system} in $\A^n_k$. Coordinate systems are transformed into each other 
by automorphisms of $\A^n_k$. 

Let $F\subset\A^n$ be locally closed, $p\in F$. 
We say that $F\subset\A^n$ is a {\em local coordinate subvariety} at $p$
if there is an open neighborhood $U\subset\A^n$ such that
$U\cap F$ is the zeroset of a subset of a coordinate system.
See also Example~\ref{ex:circ} below, where it is shown that the
circle $x^2+y^2-1=z=0$ in $\A^3$ is a local coordinate subvariety
at any of its points.

\begin{proposition} \label{thm:emb} 
Any smooth variety $Z$ admits a covering by open subsets isomorphic 
to smooth, locally closed hypersurfaces in affine space. 
If $Z$ is locally closed in $\A^n$  
and $q\in Z$, then there exists a local coordinate subvariety $F$ 
at $q$ containing $Z$ as a hypersurface.
\end{proposition}


\begin{proof}
Let $Z\subset U\subset\A^n$, $q\in Z$, and define $m:=\dim(Z)+1$.
Let $L$ be a generic element in the Grassmann variety $G(n-m,n)$,
and let $\pi:\A^n\to\A^m$ be a linear projection map that maps
$L$ to a point (e.g., the canonical map to the quotient space
$\A^n/L\cong \A^m$). Then the Zariski closure $H$ of $\pi(Z)\subset\A^m$ 
is a hypersurface. We claim that $\pi|_Z:Z\to H$ is a local
isomorphism at $q$. This is a strengthening of the theorem
of the primitive element, which says that $Z$ is birationally
equivalent to a hypersurface. (The strengthening is that we have a birational
equivalence that is \'etale at a specific point.)

To prove the claim, we form the cone of $Z$ at $q$, which is the 
Zariski closure of the union of $Z$-secants through $q$ in $\A^n$.
Its dimension is at most $m$, hence its intersection with $L$ 
is just a single point (namely $q$). 
Let $p:=\pi(q)$. 
As $L$ is generic and $k$ is infinite, the scheme-theoretic fiber 
$(\pi|_Z)^{-1}(p)$ is just $q$, and this holds also geometrically, i.e., the
residue field of $q$ is equal to the residue field of $p$.
Generic projections to subsets preserving the dimension are finite maps
(compare with the proof in \cite{Atiyah_Macdonald:69}, p. 69, of the
N\"other normalization lemma). Hence the map
$\pi|_Z$ corresponds to a local homomorphism
$i:{\mathcal O}_{H,p}\hookrightarrow{\mathcal O}_{Z,q}$ such that
${\mathcal O}_{Z,q}$ is a finite ${\mathcal O}_{H,p}$-module. Because
$\pi|_Z$ is \'etale at $q$, the completion
$\widehat{\mathcal O}_{H,p}\hookrightarrow\widehat{\mathcal O}_{Z,q}$
is an isomorphism. But completion is a faithfully exact functor on
finitely generated ${\mathcal O}_{H,p}$-modules, hence $i$ is an isomorphism.
It follows that $\pi|_Z$ is an isomorphism on suitably chosen
open neighborhoods $V$ of $q$ and $W$ of $p$.

Strictly speaking, the above argument shows only that $\pi|_Z$ gives
an isomorphism after extending the base field to the residue field of $q$,
because the construction of the cone at $q$ already requires this extension.
But since the map $\pi|_Z$ is an isomorphism after base field extension,
it must have been an isomorphism also before base field extension.

The map $\alpha:=(\pi|_Z)^{-1}$ is defined in $W\subset H$
and can be extended to an open neighborhood $W'$ of $p$ in $\A^m$. The map 
\[ \phi:W'\times L\to\pi^{-1}(W'), \ (u,l)\mapsto\alpha(u)+l \]
is then an isomorphism mapping the coordinate
subvariety $W'\times\{0\}$ to a subvariety $F$ of dimension $m$
containing $Z$.
\end{proof}

\begin{example}
Consider the space curve $Z$ given as the complete intersection 
of the two surfaces $S: y^2 = x^3-x$ and $T: z^2 = y^3 -y$. 
Both $S$ and $T$ are cylinders over the elliptic curve and therefore smooth. 
The intersection is transversal so that $Z$ itself is smooth. 
Let $U$ be the complement of $V_{\A^3} (x^3-x-1)$ in $\A^3$. 
This is an open neighborhood of the origin $p=0$. We claim that inside $U$, 
the curve $Z$ is isomorphic to a plane curve $C$. 
Indeed, substituting in the second equation $y^2$ by $x^3-x$ yields on 
$Z\cap U$ the relation $y={z^2 \over x^3-x-1}$.  
This equation defines a smooth surface (a graph) 
in $U\subseteq \A^3$ isomorphic to the open affine chart 
$\A^2 \setminus V_{\A^2} (x^3-x-1)$. 
And in this way $Z$ is isomorphic to its image curve there.
\end{example}


\section{Plain Varieties}
\label{sec:plain}

Let $V$ be a plain variety. The open subsets isomorphic to open subsets of affine space will be called {\em plain charts} of $V$.  

Any smooth toric variety is plain: A
toric variety is covered by affine toric varieties,
and a smooth affine toric variety is the product
of an affine space and a torus (see \cite{Fulton:93}, p. 29). 
Since the torus can be embedded into affine space as an open subset,
the assertion follows.

The notion of plain varieties is clearly a local one;
here is an algebraic criterion for being plain.


\begin{proposition}
A variety is plain if and only if the stalks of its structure sheaf 
are $k$-isomorphic to localizations of a polynomial algebra 
$k[x_1,\dots,x_n]$ at prime ideals.
\end{proposition}


\begin{proof}
The ``only if'' direction is obvious. Conversely, assume that $X$ is a variety,
$p$ a point in $ X$, and that the local ring ${\mathcal O}_{X,p}$ is isomorphic 
to ${\mathcal O}_{\mathbb{A}^n,q}$ for some $q\in\mathbb{A}^n$. Then the 
quotient fields of these two local rings are isomorphic over  $k$. 
Hence there exist rational maps  $\phi:X\to \mathbb{A}^n$ and 
$\psi:\mathbb{A}^n\to X$, inverse to each other, and defined at $p$ and $q$ 
respectively; and open neighborhoods $U$ of $p$ and $V$ of $q$ 
between which $\phi$ and $\psi$ induce $k$-isomorphisms. 
\end{proof}


\begin{proposition}
  Assume that $k$ is algebraically closed. 
  Then every rational smooth curve or surface is plain.
\end{proposition}


\begin{proof}
Assume that the curve $C$ is rational and smooth. Then $C$ is an
open subvariety of a complete rational smooth curve, which is
obtained by projectivization and subsequent desingularization.
As every birational map of complete smooth curves is an isomorphism 
(see \cite{Hartshorne:77}, II.6.7, p. 136), $C$ is an open subvariety 
of $\mathbb{P}^1$. (See also \cite{Hartshorne:77}, Ex. I.6.1, p. 46, 
Cor. I.6.10, p. 45, and IV.1.3.5, p.  297.)

Assume that the surface $S$ is rational and smooth. 
As in the curve case, it suffices to prove the statement for complete surfaces 
(compare with \cite{Hartshorne:77}, Rem.~II.4.10.2, p.~105).
It is well-known that any complete rational surface can be obtained
by repeated point blowups of a minimal rational surface.
The minimal rational surfaces are $\mathbb{P}^2$ and the Hirzebruch
surfaces $F_n$, $n=0,2,3,\dots$; these are all plain varieties.
The blowup of a plain surface at a point is also plain, because
the blowup of $\A^2$ at the origin can be covered by two open sets
which are both isomorphic to $\A^2$. 
\end{proof}


Obviously, plain varieties are smooth and rational. We conjecture
that also the converse holds: Any smooth and rational variety is plain.

In the following section  we will show that the class of plain varieties
is closed under blowups along nonsingular subvarieties. In order to
prove the conjecture, it would suffice to show that plain varieties are also
stable under the inverses of such blowups (blowdowns), at least in characteristic zero.
The reason is that any birational map is a composition of such blowups
and their inverses, by the Weak Factorization Theorem \cite{Abramovich_et_al:02,Wlodarczyk:03}.
The stability under blowdowns would follow if one could show that plain in codimension~1 (i.e., the non-plain
locus has codimension at least~2) implies plain.

\begin{example}
The surface $S$ with defining equation $x-(x^2+z^2)y$ over $\Q$ 
is rational---the equation is linear in $y$---and nonsingular at the origin. 
By the conjecture, there should be a plain chart containing the origin. 

Such a chart does indeed exist, namely the intersection
of $S$ with the complement of $V_{\A^3} (1-xy)$. It is isomorphic
to the complement of $V_{\A^2} (u^2v^2+1)$, by the isomorphism
\[ (u,v) \mapsto ( u^2v/(u^2v^2+1), v, u/(u^2v^2+1) ) \]
with inverse
\[ (x,y,z) \mapsto ( z/(1-xy), y ) . \]
\end{example}


\section{Blowups}

We recall the definition of blowups 
(following \cite{Hartshorne:77} and \cite{Eisenbud_Harris:01}):
Let $X$ be a variety, and let ${\mathcal I}\subseteq {\mathcal O}_X$ be a
sheaf of ideals on $X$. Then the blowup $\widetilde{X}$ of $X$ 
with center ${\mathcal I}$
is defined as the $\mathrm{Proj}$ of the sheaf of graded algebras
$ {\mathcal S} = \oplus_{i=0}^{\infty}\, {\mathcal I}^i$. 
If $X$ is affine, say $X=\mathrm{Spec}(R)$, and 
$\mathcal I$ is the sheaf corresponding to the ideal $I\subseteq R$
that is generated by $a_1,\dots,a_m\in R$, then $\widetilde{X}$ can be 
embedded into $\P_R^m$ as the closed subvariety defined by the
homogeneous ideal $J=\ker(\varphi)$, where 
$\varphi:R[y_1,\dots,y_m]\to R[t]$ is defined as the $R$-homomorphism
sending $y_i$ to $ a_it$ for $i=1,\dots,m$.
Here is a description of this situation in terms of affine charts.


\begin{lemma}\label{lem:regular}
Let $X,I$ and $a_1,\dots,a_m$ be as above.
Then the blowup $\widetilde{X}$ is covered by affine open subsets $U_1,\dots, U_m$ with $U_i$ isomorphic to
$\mathrm{Spec}(R[\frac{a_1}{a_i},\dots,\frac{a_m}{a_i}])$.
\end{lemma}


\begin{proof}
For $i=1,\dots,m$, we define $U_i$ as the intersection of $\widetilde{X}$
with the affine patch $y_i\ne 0$ in $\P_R^m$. It has affine coordinates
$$x_1=\frac{y_1}{y_i},\dots,x_{i-1}=\frac{y_{i-1}}{y_i},x_{i+1}=\frac{y_{i+1}}{y_i},\dots,x_m=\frac{y_m}{y_i}$$
(in addition to the affine coordinates for $X$).
Then the defining ideal of $U_i\subset \A_R^m$ is the image of $J$
defined as above under the dehomogenization homomorphism
$R[y_1,\dots,x_m]\to R[x_1,\dots,x_{i-1},x_{i+1},\dots,x_m]$ mapping
$y_i$ to $1$ and $y_j$ to $x_j$ for \linebreak $j\ne{i}$. But this ideal is also
the kernel of the ring homomorphism \linebreak
$R[x_1,\dots,x_{i-1},x_{i+1},\dots,x_m]\to \rm {Quot}(R)$  mapping $x_j$ to $a_j/a_i$,
hence the quotient ring is isomorphic to 
$R[\frac{a_1}{a_i},\dots,\frac{a_m}{a_i}]$.
\end{proof}


Assume that $X={\mathbb A}^n_k$  and that $I$ is the ideal of a coordinate
subvariety, e.g., $I=(x_1,\dots,x_r)$.
Then, for $i=1,\dots,r$, the ring $R[\frac{x_1}{x_i},\dots,\frac{x_r}{x_i}]$ 
is isomorphic to the polynomial ring $k[y_1,\dots,y_n]$ by the isomorphism 
sending $\frac{x_j}{x_i}$ to $y_j$ for $j\ne i$ and $x_i$ to $y_i$. 
Hence the spectra $U_i$ of these rings are affine spaces, and $\widetilde{X}$ is plain. As being plain is a local property, this shows that the blowup of $\A^n_k$ along local coordinate subvarieties is plain.


\begin{example} \label{ex:circ}
The circle given by the ideal $I=(x^2+y^2-1,z)$ in $\A^3$ is a local
coordinate subvariety. For instance, let $q:=(1,0,0)$. Then the rational map
\[ (x,y,z) \mapsto (u,v,w)=\left(\frac{x}{y+1},
	\frac{x^2+y^2-1}{y+1},z\right) \]
is an isomorphism of an open neighborhood of $q$ to an open neighborhood of $(0,0,0)$; its birational
inverse is
\[ (u,v,w) \mapsto (x,y,z) = \left(\frac{2u}{u^2-v+1},
	\frac{-u^2+v+1}{u^2-v+1},w\right) . \]
Any other point on the circle can be moved to $q$ by a rotation.
This shows directly (i.e., without referring to the theorem below) that the blowup of the circle is plain (see also 
\cite{Bodnar_Schicho:98}).
\end{example}

In order to prove the stability of plain varieties under blowup in general, the observation above is not enough,
for there are smooth subvarieties in $\A^n_k$ which are not
local coordinate varieties, e.g., non-rational curves.


\begin{theorem} \label{thm:algebraic}
Let $F$ be a coordinate subvariety of $U$ open in $\A^n_k$,
and let $Z$ be a closed smooth hypersurface of $F$. Assume that $k$ is infinite.
Then the blowup of $U$ along $Z$ is plain.
\end{theorem}


\begin{proof}
Assume that $F$ is defined by $x_1=\dots=x_r=0$, and that $Z$ is defined
within $F$ by $f(x_{r+1},\dots,x_n)=0$. Let $p\in Z$ be a point,
and assume that $\partial_{x_n}{f}$ does not vanish
at $p$. We will construct an open neighborhood $V$ of $p$ in $U$ and a covering
of the blowup $\widetilde V$ of $V$ along $Z\cap V$ by plain charts.

Note first that 
$\mathrm{Spec}(k[U][\frac{x_1}{f},\dots,\frac{x_r}{f}])$ is a plain chart
of the blowup. Indeed, the ring 
$k[x_1,\dots,x_n,\frac{x_1}{f},\dots,\frac{x_r}{f}]$ is isomorphic to
$k[y_1,\dots,y_r,x_{r+1},\dots,x_n]$ by the isomorphism fixing
$x_{r+1},\dots,x_n$ and sending $x_i$ to $y_if(x_{r+1},\dots,x_n)$
for $i=1,\dots,r$. Unfortunately, the other charts of the blowup as described in
Lemma~\ref{lem:regular} need not be plain charts in general. The trick is now to choose other defining equations for $Z$, thus changing the charts on $\widetilde V$.

For $i=1,\dots,r$, we set $f_i=f(x_{r+1},\dots,x_{n-1},x_n+x_i)$.
We claim that there is an open neighborhood $V$ of $p$ 
where $f,f_1,\dots,f_r$ are generators of the ideal of $Z$.
This is a consequence of the Taylor expansion
\[ f_i = f+x_i\cdot\left(\partial_{x_n}{f}+h_i\right) \]
for suitable polynomials $h_i$ vanishing at $p$ ($i\leq r$);
we choose the affine open neighborhood $V$ of $p$ so that the polynomials 
$g_i:=\partial_{x_n}{f}+h_i$ do not vanish in $V$. They 
therefore have multiplicative inverses in $k[V]$.

The first  chart of $\widetilde V$ is isomorphic to the spectrum of $k[V][{f_1\over f},\ldots,{f_r\over f}]$. From ${f_{i}\over f}= 1+g_i\cdot {x_i\over f}$ for $i$ from $1$ to $r$ and because all $g_i$ are invertible on $V$ it follows that $k[V][{f_1\over f},\ldots,{f_r\over f}]=k[V][{x_1\over f},\ldots,{x_r\over f}]$. This algebra is a polynomial ring and therefore the chart is plain. 

 Any of the other charts of $\widetilde V$ is isomorphic to the spectrum of $k[V][{f\over f_i},{f_1\over f_i},\ldots,{f_r\over f_i}]$ for some $i$ between $1$ and $r$. This algebra equals $k[V][{x_1\over f_i},\ldots,{x_r\over f_i}]$, because of the relations ${f\over f_i}= 1-g_i\cdot {x_i\over f_i}$ and ${f_j\over f_{i}}= {f\over f_i}+ g_j\cdot {x_j\over f_i}$, respectively ${x_j\over f_i}=g_j^{-1}\cdot ({f_j\over f_{i}}- {f\over f_i})$.
 Hence these other charts are also plain. The theorem is proven.
\end{proof}


\noindent As a consequence of Theorem~\ref{thm:algebraic} we obtain the stability of plain varieties under blowup.

\begin{theorem} \label{thm:main}
Over infinite fields, the blowup of a plain variety along a smooth subvariety is plain.
\end{theorem}

\begin{proof} 
As being plain is a local property, and using the local embedding
of smooth varie\-ties from Proposition \ref{thm:emb}, it suffices to consider blowups along smooth hypersurfaces
in coordinate subvarieties; this is just the case settled in
Theorem~\ref{thm:algebraic}.
\end{proof}


\begin{example} 
We blow up $\A^3$ with center the plane elliptic curve $Z:z=x-x^3+y^2=0$. We place ourselves in  some
neighborhood $V$ of the origin of $\A^3$. 
Setting $f_2=x-x^3+y^2$, $f_1=x+z-(x+z)^3+y^2$, and $g=1-3x^2-3xz-z^2$,
we get $f_1=f_2+gz$. Hence $f_1$ and $f_2$ generate the ideal of $Z$
in the principal open set $V$ defined by $g\ne 0$.

In the first chart $V_1$, we adjoin the fraction $v=f_2/f_1$ to $k[V]$.
We introduce the new variable $w$ for $x+z$ and use it for eliminating $x$;
in particular, we write $g$ in the form $g(w,y,z)=1-3w^2+3wz-z^2$.
Let the variable $s$ stand for the fraction $z/f_1=g^{-1}(1-v)$.
Express now $z$ as $tf_1(w,y)$ and $v$ as
$1-tg(w,y,z)=1-tg(w,y,tf_1(w,y))$. This chart is isomorphic to the
open set in $\A^3$ with coordinates $w,y,t$ defined by the inequality $g(w,y,tf_1(w,y)) \ne 0$, say
\[ 1-3w^2+3wt(w-w^3+y^2)-t^2(w-w^3+y^2)^2 \ne 0 , \]
and the exceptional divisor inside this chart has equation
\[ f_1(w,y) = w-w^3+y^2 =0. \]

\noindent In the second chart $V_2$ of the blowup, we adjoin the
fraction $u:=f_1/f_2$ to $k[V]$. We introduce a new variable $s$ for the
fraction $z/f_1=g^{-1}(u-1)$. Now, we can express $z$ as $sf_2(x,y)$
and $u$ as $1+sg(x,y,z)=1+sg(x,y,sf_2(x,y))$. The chart is
isomorphic to the open set of $\A^3$ with coordinates $x,y,s$
defined by the inequality $g(x,y,sf_2(x,y)) \ne 0$, say
\[ 1-3x^2-3xs(x-x^3+y^2)-s^2(x-x^3+y^2)^2 \ne 0 , \]
and the exceptional divisor inside this chart has equation
\[ f_2(x,y) = x-x^3+y^2 =0. \]

\noindent The intersection $V_1\cap V_2$ is isomorphic to the open subset 
$V_2'\subset V_2$ defined by $u=1+sg(x,y,sf_2(x,y))  \ne 0$, say
\[ 1+s-3x^2s-3xs^2(x-x^3+y^2)-s^3(x-x^3+y^2)^2 \ne 0 \]
and to the open subset $V_1'\subset V_1$ defined by $v=1-tg(w,y,tf_1(w,y)) \ne 0$, say
\[ 1-t+3w^2t-3wt^2(w-w^3+y^2)+t^3(w-w^3+y^2)^2 \ne 0 . \]
The chart change map $V_2'\to V_1'$ is given by
\[ (x,y,s) \mapsto (w,y,t)=
	\left(x+sf_2(x,y),y,\frac{s}{1+sg(x,y,sf_2(x,y))}\right) = \]
\[ \left(x+xs-x^3s+y^2s,y,\frac{s}{1+s-3x^2s-3xs^2(x-x^3+y^2)-s^3(x-x^3+y^2)^2}\right) . \]

\[\]
\end{example}



\end{document}